\documentclass{article}
\usepackage{graphicx}
 \usepackage{mathptmx}
\usepackage{amsmath, amstext,amssymb,amsfonts}
\usepackage[english]{babel}
\usepackage{delarray}
\usepackage{mathptmx}
\usepackage{amsmath, amstext,amssymb,amsfonts}
\usepackage[english]{babel}
\usepackage{delarray}

\DeclareFontFamily{U}{mathx}{\hyphenchar\font45}
\DeclareFontShape{U}{mathx}{m}{n}{
      <5> <6> <7> <8> <9> <10>
      <10.95> <12> <14.4> <17.28> <20.74> <24.88>
      mathx10
      }{}
\DeclareSymbolFont{mathx}{U}{mathx}{m}{n}
\DeclareFontSubstitution{U}{mathx}{m}{n}
\DeclareMathAccent{\widecheck}{0}{mathx}{"71}
\DeclareMathAccent{\wideparen}{0}{mathx}{"75}

\numberwithin{equation}{section}

\newtheorem{thm}{Theorem}[section]

 \newtheorem{prop}[thm]{Proposition}
 
  \newcommand{\R}{\mathbb R}    \newcommand{\Co}{\mathbb C}             \newcommand{\Z}{\mathbb Z}  \newcommand{\M}{M({\mathbb R})}
        \newcommand{\CH}{CH(\mathbb R)}      \newcommand{\Le}{L^1({\mathbb R})}     \newcommand{\cn}{c_n}
\date{}
\title{A note on powers of the characteristic  function}
\author{{\large Saulius  Norvidas} }
\date{{\footnotesize Institute of Data Science and Digital Technologies, Vilnius University, Akademijos str. 4, Vilnius LT-04812, Lithuania\\
 ({\rm{e-mail: norvidas{@}gmail.com}})}}
\begin{document}

\maketitle
 {{ {\bf Abstract}}}
 Let $\CH$ denote the family of  characteristic functions of probability measures (distributions) on the real line $\R$. We study the following  question: given an integer $n>1$,  do  there  exist  two  different  $f, g\in\CH$ such  that $ f^n\equiv g^n$? For positive even  $n$, well-known examples answer this question in the affirmative. It turns out that the same  is true also for any odd $n>1$. For $f\in CH(\R)$ and integer $n>1$, set $C_n(f)=\{g\in CH(\R): g^n\equiv f^n\}$. In  this paper,   we give an estimate for cardinality  (or cardinal number) of  $C_n(f)$. In addition, we describe such $f$ for which our estimate is sharp.



\section{Introduction}\label{s:1}
{\large{
 Let  $M(\R)$ denote  the  Banach algebra of bounded regular   Borel measures  on the real line $\R$  with the total variation norm $\|\mu\|$ and convolution as multiplication.  We define the Fourier-Stieltjes transform of $\mu\in \M$   by
\[
\hat{\mu}(x)=\int_{\R}e^{-it x}\, d\mu(t),
\]
$x\in\R$.

If $\mu\in \M$ is  positive on the whole $\R$ and  $\|\mu\|=1$, then  in the language of probability theory,  such an $\mu$ is called a  probability measure.   If $\mu$ is a probability measure,  then $f(x):= \hat{\mu}(-x)$, $x\in\R$,  is  called the  characteristic function of $\mu$. In the sequel,  $\CH$ denotes the multiplicative semigroup of characteristic functions on $\R$.

 Let   $\Z$  denote  the usual  group of integers. Given $f\in \CH$ and $n\in\Z$, $n\ge 2$, we define
\[
C_n(f)=\{g\in\CH: g^n\equiv f^n\}.
\]
 We say that $\cn(f)$ is  trivial if   $\cn(f)$ do not contains elements of $\CH$ other than $f$. The purpose  of this paper is to study the structure of $\cn(f)$. Our  interest in $\cn(f)$ is inspired by the following question in \cite[p.  334]{7} as an unsolved problem:
\begin{gather}
Do \ there \ exist \ two \ different \ characteristic \ functions \ f \ and \ g \ such \ that \nonumber\\
  f^n\equiv g^n \ for \ some \ odd \ n, \ n>1?
\end{gather}
If there $n\in\Z$ is   positive and even, then the answer to this question is well known. Indeed, let  $f$ and $g$ be periodic functions with periods $2$ and $4$, respectively, such that $f(x)=1-|x|$, for $|x|\le 1$, and $g(x)=1-|x|$, for $|x|\le 2$. Then $f$ and $g$ are characteristic functions (see, e.g., \cite[p. 265]{7}) such that $|f|\equiv |g|$. Therefore, $f^{2k}\equiv g^{2k}$ for all positive $k\in\Z$.

 If  $n>1$ and $n$ is odd, then an affirmative answer to (1.1) comes from the following theorem.
\begin{thm} (see \cite[theorem 1.1]{5}). There exists  $f\in CH(\R)$ such that $C_n(f)$ is non-trivial for all integer $n>1$.
\end{thm}
In this paper, we will study the structure of  $C_n(f)$ in more detail. It turns out that in such a problem, the  significant role plays the (geometrical)  construction  of the null set $\{x\in\R: f(x)=0\}$ of $f$.

For $f\in \CH$, let us define the essential support of $f$ by $S_f=\{ x\in\R: \ f(x)\neq 0\}$. Obviously,  $S_f$ is an open subset of $\R$. It is known that  for any open subset $U\subset \R$ such that $-U=U$ and $0\in U$, there exists $f\in \CH$ with $S_f=U$ (see, e.g., \cite{2}).  Next, $S_f$   can  be represented as a finite or infinite  union $\bigcup_{j}E_j$, where $\{E_j\}_j$ is the family of all open connected components (intervals) of $S_f$. The cardinal (or cardinal number) ${\text{card}}(\{E_j\}_j)$ of $\{E_j\}_j$ we denote by ${\text{comp}}\,(S_f)$. Of course,  if $\{E_j\}_j$ is an infinite family, then ${\text{comp}}\,(S_f)=\infty$. Since  any $f\in CH$  satisfies
\begin{equation}
f(-x)=\overline{f(x)}
\end{equation}
for all  $x\in\R$, it follows that if  ${\text{comp}}\,(S_f)<\infty$, then there is $k\in\Z$, $k\ge 0$, such that  ${\text{comp}}\,(S_f)=2k+1$. Moreover,  then  there exists a finite sequence of pairwise disjoint open  intervals $E_0=(-b_0,b_0)$ and $E_j=(a_j,b_j)$,  $\ E_j=-E_{-j}$, $j=-k,\dots,-1,1,\dots, k$,  such that
\begin{equation}
S_f=\bigcup_{j=-k}^{k}E_j
\end{equation}
with
\begin{equation}
 \ b_0<a_1<b_1<a_2<b_2<\dots <a_k<b_k .
\end{equation}
Note that  in (1.3)  also might be $b_k=+\infty$.

\begin{prop}
 Let $f\in \CH$ and let $n\in\Z$, $n\ge 2$. Suppose that  ${\text{comp}}\,(S_f)=2k+1$  with certain $k\in\Z$, $k\ge 0$. Then
\begin{equation}
{\rm{card}}(\cn(f))\le n^k.
\end{equation}
\end{prop}

\begin{thm}
 Under the conditions of Proposition 1.1,  for each $n\in\Z$, $n\ge 2$ and any sequence of of pairwise disjoint open  intervals  $\{E_j\}$   satisfying (1.3) and (1.4),   there exists $f\in \CH$ such that\begin{equation}
{\rm{card}}(\cn(f))= n^k.
\end{equation}
\end{thm}
{\bf{Remark 1.3}}. \ Note a particular case of Proposition 1.1 when ${\text{comp}}\,(S_f)=1$, i.e., if $k=0$. For example, the characteristic functions of the following  frequently used probability measures (distributions) satisfy this condition ($k=0$): Normal, Laplace, Poisson, Cauchy and some other (see, for example, \cite[p.p. 282-329]{7}). Therefore, (1.5) shows  that these $f$   generate the trivial classes $C_n(f)$  for all integer $n > 1$.

\section{ Preliminaries and proofs }
The Lebesgue space $\Le$  can be identified with the closed ideal in $\M$ of   measures absolutely continuous with respect to the Lebesgue measure $dx$ on $\R$. Namely, if $\varphi\in \Le$, then  $\varphi$ is associated with  the  measure
\[
\mu_{\varphi}(E)=\int_{E}\varphi(t)\,dt
\]
for  each Borel subset $E$ of $\R$. Hence  $\widehat{\varphi}(x)=\int_{\R}e^{-itx}\varphi(t)\,dt$.  In particular, if $\mu=\varphi(t)dt$, where  $\varphi\in L^1(\R)$ and $\varphi$ is such  that  $\|\varphi\|_{L^1(\R)}=1$ and $\varphi\ge 0$ on $\R$,   then $\varphi$  is called the probability density function of $\mu$, or the probability density for short.

We  define the inverse Fourier transform by
 \[
\widecheck{\psi}(t)=\frac1{2\pi}\int_{\R}e^{it\,x}\psi(x)\,dx,
\]
$t\in\R$. Then  the inversion  formula  $\widehat{(\check{\psi})}=\psi$ holds for suitable $\psi\in L^1(\R)$.  A function $F:\R\to\Co$ is said to be positive definite if
\[
\sum_{i,j=1}^m F(x_j-x_i)c_j\overline{c_i}\ge 0
\]
holds for all finite sets  $c_1,\dots,c_m\in\Co$ and $x_1,\dots,x_m\in\R$. The Bochner theorem (see, e.g., \cite[p. 71]{4}) states that a continuous function $F:\R\to\Co$ is positive definite if and only if there exits a positive measure $\mu_F\in M(\R)$ such that $F=\widehat{\mu_F}$. Hence, any $f\in \CH$ is positive definite with $f(0)=1$. In the case  if $\varphi:\R\to\Co$ is continuous and $\varphi\in L^1(\R)$, then $\varphi$ is positive definite if and only if $\widecheck{\varphi}(t)\ge 0$ for all $t\in\R$.

{\bf{Proof of Proposition 1.1.}}\quad  Let $g\in\cn(f)$. Then it is immediate that $S_f=S_g$ and
\begin{equation}
|f(x)|=|g(x)|,
\end{equation}
$x\in S_f$. Fix any $E_j\subset S_f$ in (1.3). Using the fact that $E_j$ is an open connected component of $S_f$, we have that there are two continuous functions $u_{f,j}, u_{g,j}: E_j\to \R$ such that
\begin{equation}
f(x)=|f(x)|e^{iu_f(x)}\quad {\text{and}}\quad  g(x)=|g(x)|e^{iu_g(x)}
\end{equation}
for all  $x\in E_j$. Since $f^n\equiv g^n$ and $u_{f,j}$ together with $u_{g,j}$ are continuous on  $E_j$, it follows from (2.1) and (2.2) that
\begin{equation}
u_{g,j}(x)=u_{f,j}(x)+2\pi m_j
\end{equation}
for certain $m_j\in\{0,1,\dots,n-1\}$ and all $x\in E_j$. Therefore,
\begin{equation}
g(x)=\sum _{j=-k}^{k}f(x)\chi_{E_{j}}(x) e^{i2\pi m_j},
\end{equation}
$x\in \R$, where   $\chi_{E_{j}}$ is the indicator function of $E_j$.    (1.2) implies  that $m_{-j}=-m_j$ for each $j=0,1,\dots, k$. Next, we claim that $m_j=0$. Indeed, this claim follows from (2.3) and equalities $u_{f,j}(0)=u_{g,j}(0)=0$, since $f(0)=g(0)=1$.  Hence, we conclude from (2.4) that
\begin{equation}
g(x)=f(x)\chi_{E_0}(x)+\sum_{j=1}^{k}f(x)\Bigl(\chi_{E_j(-x)}(x)e^{-i2\pi m_j}+\chi_{E_j(x)}(x)e^{i2\pi m_j}\Bigr)
\end{equation}
for all $x\in \R$. Keeping in mind that each $m_j$, $j=1,2,\dots,n-1$, may take any  value in $\{0,1,\dots,n-1\}$, we obtain from (2.5) the estimate (1.5).  Proposition 1.1 is proved.

{\bf{Remark 2.1}}. \
In Remark 1.3 and in the proof of Proposition 1.1 (see (2.1)) we concerned with certain aspects of the phase retrieval problem, i.e., with the problem of recovery of a measure $\mu$ (not necessarily a probability measure)   when are given the amplitude $|\widehat{\mu}|$ of its Fourier transform. This problem is well known in various fields of science (see survey \cite{6}). The question about uniqueness of the solution for such problem arises in quantum mechanics (see  \cite{3}).

{\bf{Proof of Theorem 1.2.}}\quad

Let us split the proof into two cases. First  we consider the case when all $E_j$ in (1.2) are finite intervals. Let us denote by $|E_j|$ and by $2\varrho$  the length of interval $E_j$ in (1.3) - (1.4) and the minimal of these lengths, respectively.  Since $a_0=-b_0$, i.e., $E_0=(-b_0,b_0)$, it follows from (1.4) that
\begin{equation}
2\varrho=\min\{ 2b_0; \ \min_{1\le j\le k}(b_j-a_j)\}.
\end{equation}
 Now fix any $\varphi\in \CH$ such that
\begin{equation}
{\text{supp}}\,\varphi\subset [-\varrho,\varrho] \quad \ {\text{and}} \quad \ \varphi(x)>0 \ \ {\text{for \ all}} \ x\in (-\varrho,\varrho).
\end{equation}
Since, in addition,  $\varphi\in L^1(\R)\cap L^{\infty}(\R)$, we have  that $\widecheck{\varphi}\in L^1(\R)$  (see, e.g., \cite[p. 409]{1}). Also the second condition in (2.7) shows  that  $\varphi$             is an even function.

 $\widecheck{\varphi}\in L^1(\R)$ with
Now fix any $j\in\{0,1,\dots,k\}$. It follows immediately from (2.6) and (2.7) that there exists a positive number $m=m(j)$ and a family of points $\{\tau_{j_{p}}\}_{p=1}^{m(j)}$    satisfying
\begin{equation}
\tau_{j_1}=a_j+\varrho<\tau_{j_2}<\dots <\tau_{j_{m}}=b_j-\varrho
\end{equation}
 such that the function
\begin{equation}
F_j(x)=\sum_{p=1}^{m(j)}\varphi(x-\tau_{j_p})
\end{equation}
is supported on $[a_j,b_j]=\overline{E_j}$. Moreover, $F_j$   is strictly positive on $(a_j,b_j)=E_j$. Since $\varphi$ is even, it follows that the  function
 \begin{equation}
 F_j(-x)=\sum_{p=1}^{m(j)}\varphi(x+\tau_{j_p})
\end{equation}
 is supported on $\overline{E_{-j}}=-\overline{E_j}$ and strictly positive on $-E_j$.

Now  take an  arbitrary sequence  of positive numbers   $\{\alpha_j\}_{j=0}^k$ such that
\begin{equation}
\alpha_j\le \frac1{2m_j(k+1)},
\end{equation}
$j=0,1,\dots,k$,  and define
\begin{equation}
F(x)=\varphi(x)+\sum_{j=0}^k\alpha_j\Bigl(F_j(x)+F_j(-x)\Bigr).
\end{equation}
According to (2.6), (2.7), (2.8) and properties of the  functions $F_j$ mentioned above, we see that
\begin{equation}
S_F=\bigcup_{j=-k}^{k}E_j .
\end{equation}
We claim that $F$ is continuous and  positive definite. Indeed, since $F$ is continuous and compactly supported, we it follows that  $F$ and $\widecheck{F}$ are in $L^1(\R)$. Then we deduce from  (2.9) and (2.10) that
\begin{gather}
\widecheck{F}(t)=\widecheck{\varphi}(t)+\sum_{j=0}^k\biggl(\alpha_j\widecheck{\varphi}(t)\sum_{p=1}^{m(j)}2\cos(\tau_{j_p}\cdot t)\biggr)
  =\widecheck{\varphi}(t)\biggl[1+2\sum_{j=0}^k\biggl(\alpha_j\sum_{p=1}^{m(j)}\cos(\tau_{j_p}\cdot t)\biggr)\biggr].
 \end{gather}
By combining  that $\varphi\in \CH$ with condition (2.11), we get that $\widecheck{F}(t)\ge 0$ for all $t\in\R$. Hence,   Bochner's theorem proving our claim.

Next, let $n\in\Z$, $ n\ge 2$. For any $\omega\in\Z^{k}$, let us define
\begin{gather}
G_{\omega}(x)=\varphi(x)+\alpha_0\Bigl(F_0(x)+F_0(-x)\Bigr)
+\sum_{j=1}^k\alpha_j\biggl[e^{i2\pi \omega_j/n}F_j(x)+e^{-i2\pi \omega_j/n}F_j(-x)\Bigr)\biggr].
\end{gather}
Then analogous to (2.10), we obtain
\begin{equation}
\widecheck{G_{\omega}}(t)=\widecheck{\varphi}(t)\biggl(1+2\alpha_0\sum_{p=1}^{m(0)}\cos(\tau_{0_p}\cdot t)+2\sum_{j=0}^k\biggl[\alpha_j\sum_{p=1}^{m(j)}\cos\Bigl(\frac{2\pi \omega_j}{n}+\tau_{j_p}\cdot t\Bigr)\biggr]\biggr).
\end{equation}
Therefore,  using again (2.11), we see that  each such an $G_{\omega}$ is also continuous and positive definite.

Now we claim that
\begin{equation}
G_{\omega}^n\equiv F^n
\end{equation}
for all $\omega\in\Z^{k}$. This follows immediately from definitions $F$ in (2.12) and $G_{\omega}$ in (2.15). Indeed, at first, ${\text{supp}}\,F= {\text{supp}}\,G_{\omega}=\cup_{j=-k}^{k}E_j$. Next,  it is easy to check that
\[
\begin{array}\{{lcr}.
G_{\omega}(x)=\alpha_je^{i2\pi \omega_j/n}F_j(x)=e^{i2\pi \omega_j/n}F(x), \qquad \ \ \ {\text{for }} \   x\in E_j, j=1,\dots, k, \\
G_{\omega}(x)=\varphi(x)+\alpha_0(F_0(x)+F_0(-x))=F(x), \quad \ {\text{for }} \   x\in E_0,\\
G_{\omega}(x)=\alpha_je^{-i2\pi \omega_j/n}F_j(-x)=e^{-i2\pi \omega_j/n}F(x),  \quad {\text{for }} \   x\in E_j, j=-1,\dots,- k.
\end{array}
\]
In particular, from here  and from (2.15) we obtain
\[
G_{\omega}(0)=F(0)=1+2\alpha_0\sum_{p=1}^{m(0)}\varphi(\tau_{0_p})>1.
\]
Set
\[
f(x)=\frac{F(x)}{F(0)} \quad {\text{and}} \quad g_{\omega}(x)=\frac{G_\omega(x)}{G_\omega(0)}.
\]
Then $f, g_\omega\in CH$ and $g_\omega\in C_n(f)$ for all $\omega\in\Z^k$. Moreover, (2.13) implies that $S_f=S_F=\cup_{j=-k}^{k}E_j$. Finally, it is easy seen that $G_{\omega}=G_{\theta}$ for some $\omega$ and $ \theta$ from $\Z^k$,  if and only if $\omega-\theta\in n\Z^k$. Therefore, when $\omega$ runs in $\Z^k$, then the family (2.14)  and corresponding family $\{g_\omega\}_{\omega\in\Z^k}$ contains exactly $n^k$ different characteristic functions. This shows (1.6).

Now we consider the second case when $E_k=(a_k,\infty)$. Using the same $\varphi\in CH$ satisfying (2.7), let us define by (2.9) the functions $F_j$ for $j=0,1,\dots, k-1$. Next, take an arbitrary sequence of  positive numbers $\{\gamma_l\}_{l=1}^{\infty}$ such that
\[
\sum_{l=1}^{\infty} \gamma_l=1
\]
and define
\[
F_k(x)=\sum_{l=1}^{\infty} \gamma_l\varphi(x-a_k-l\varrho).
\]
Now for any sequence of positive numbers $\{\alpha_j\}_{j=0}^k$ satisfying (2.11) for all $j=0,1,\dots, k-1$ and such that
\[
\alpha_k\le \frac1{2(k+1)},
\]
we have  that $F_k$  is supported on $\overline{E_k}$ and is strictly positive on $E_k$. The following proof in this case  exactly the same as that of  the first case.

 }}

\begin{thebibliography}{}

\bibitem{1}Hewitt, E.,  Stromberg, K.:  Real and Abstract Analysis. A modern treatment of the theory of functions of a real variable.  Springer-Verlag, New York-Heidelberg (1975).
\bibitem{2} Ilinski\u{i}, A. I.: The zeros and the argument of a characterisitic function. (Russian. English summary). Teor. Verojatnost. i Primenen.  20, no. 2, 421-427 (1975).
\bibitem{3} Klibanov, M.V., Kamburg, V. G.: Uniqueness of a one-dimensional phase retrieval problem. Inverse Problems.   30,  no. 7, 075004, 1-10  (2014).
\bibitem{4}  Lukacs, E.:   Characteristic Functions, 2nd edn. Hafner Publishing Co., New York (1970).
\bibitem{5} Norvidas, S.:  On powers of the characteristic function.  Mediterr. J. Math. V. 17 ,   Issue: 2,     Article no. 72, 1-12 (2020)
\bibitem{6} Shechtman, Y.,  Eldar,  Y.C.,  Cohen, O.,  Chapman, H.N.,  Miao, J.,  Segev, M.: Phase retrieval with application to optical imaging. IEEE signal processing magazine.   May,  87- 109 (2015).
\bibitem{7}  N. G.  Ushakov:    Selected Topics in Characteristic Functions.  Modern Probability and Statistics., VSP, Utrecht (1999).

\end{thebibliography}
\end{document}